\documentstyle[12pt,amssymb,latexsym,amsmath]{article}
\begin{document}

\title{A variant of the Dressing Method applied to nonintegrable multidimensional
nonlinear Partial Differential Equations}

\author{
A.I. Zenchuk\\
Center of Nonlinear Studies of L.D.Landau Institute
for Theoretical Physics  \\
(International Institute of Nonlinear Science)\\
Kosygina 2, Moscow, Russia 119334\\
E-mail: zenchuk@itp.ac.ru}

\maketitle

\begin{abstract}
We describe a variant of the dressing method  giving alternative representation
of multidimensional nonlinear PDE as a system of
 Integro-Differential Equations (IDEs) for  spectral and   dressing functions.
In particular, it becomes single
linear  Partial Differential Equation (PDE) with
potentials expressed through the field of the nonlinear PDE.
The absence of linear overdetermined system associated with nonlinear PDE  creates an
obstacle to obtain evolution of the spectral data (or dressing
functions): evolution   is defined by nonlinear IDE (or PDE in
particular case).
As an example, we consider
generalization of the dressing method applicable to  integrable (2+1)-dimensional 
$N$-wave  and Davey-Stewartson equations. Although represented  algorithm 
does not supply an 
analytic particular solutions, this approach may have a perspective development. 
\end{abstract}

\section{Introduction}
Completely integrable multidimensional Partial Differential
Equations (PDEs) represent attractive subject of intensive study during
last decades after the paper \cite{GGKM}. This popularity is due to
their remarkable mathematical properties and variety of physical
applications, which may be found in literature.
Investigation approach  considered in  this paper is, in some sence, associated
with   
so-called $S$-integrable PDEs \cite{C_int0},
i.e. nonlinear PDEs which  may be "linearized" using special 
technique, such as Inverse Spectral Transform (IST)
\cite{ZMNP,AKNSB,AC}.
It is well known that IST is not the only method to study 
$S$-integrable PDEs. One may refer to    
Sato Theory \cite{DKJM,DKJM2,SS,OSTT},
 Symmetry Approach \cite{MShS,KM},
Dressing Method \cite{ZSh1,ZSh2,ZM,BM} . The later, in turn, has several
formulations: Zakharov-Shabat method \cite{ZSh1}, 
local Riemann problem \cite{ZSh2}, nonlocal Riemann and
$\bar\partial$-problem \cite{ZM,BM,K}.  

Classical $S$-integrable systems are basically (1+1)- and
(2+1)-dimensional. Only special types of multidimensional
$S$-integrable examples are known, such as 
self-dual Yang-Mills equations \cite{YM,BPST,BZ,ADHM,DM} and the
Plebanski heavenly equation \cite{P,BK,MS}.
Recently a new type of multidimensional {\it partially} integrable
systems have been found \cite{ZS}, for which integration algorithm is based
on the integral operator with nontrivial kernel, which is a variant of the dressing method. 
This recent result encourage us to search for other improvements of the dressing method.

It is well known, that dressing method has  been originally developed
 to construct nonlinear PDEs together with their solutions. 
 Variant of the dressing method suggested  here does not allow one to find analytic solutions for nonlinear PDEs. However
 \begin{enumerate}
 \item
  it gives an alternative representation of largely arbitrary  nonlinear PDE as nonlinear system of Integro-Differential Equations (IDEs). In particular case, this system becomes single linear PDE where potentials are expressed through the spectral function from one hand and through the field of original nonlinear PDE from another hand;
 \item
 it relates  a {\it single} linear spectral evolution equation (written  for some spectral function) with largely arbitrary nonlinear PDE.
 \end{enumerate}   
This is an interesting result of the paper. However, the fact that one has single linear equation associated with given nonlinear PDE (instead of overdetermined linear system, like in $S$-integrable case) results in system of nonlinear IDEs (or PDEs) defining evolution of the dressing function, which is    disadvantage of our representation. Remember that dressing functions of $S$-integrable PDE  satisfy linear PDE. As a consequence, our (largely arbitrary) PDE may not be derived as compatibility conditions of  linear overdetermined system.      
  
 In some sence, similar purpose (but different approach)
 was sought in series of papers generalizing known 
 (2+1)- and (1+1)-dimensional completely integrable equations.
 These are   
  generalization of Kadomtsev-Petviashvili
equation (KP) using deformation of the classical Inverse Spectral
Transform (IST) \cite{BB},
generalizations of Korteweg-de Vries equation (KdV) and 
Nonlinear Shr\"odinger equation (NLS)  \cite{FA}, 
generalization of Benjamin-Omo
equation (BO) \cite{KLM}. In these papers 
evolution of spectral data  is defined by nonlinear nonlocal 
equations (spectral data  are replaced by  dressing functions in our case).

Here we start with  dressing method based on the integral
  equation in the form 
 \cite{SAF,Z2}, where we  introduce an integral operator
 with different type of kernel
 allowing us  to increase dimensionality of PDE.
 As a consequence, an arbitrary function of $x_i$ (independent variables of nonlinear PDE) appears in the dressing algorithm (see function $\hat \Phi(\lambda_1;x)$ in Sec.\ref{IDE}) enforcing us to introduce an extra constrain in the form of largely arbitrary  nonlinear IDE for $\hat \Phi(\lambda_1;x)$, see eq.(\ref{condition_f}).
 Fixing function $\hat \Phi(\lambda_1;x)$, this constrain provides possibility to write single nonlinear PDE for single field $u$ expressible in terms of the dressing 
 and spectral functions. Note that similar extra constrain has been introduced in \cite{ZS}, but arbitrary function there has quite different origin.  

Below we concentrate on multidimensional 
generalizations of dressing algorithm
for (2+1)-dimensional $N$-wave equation  and Davey-Stewartson equation (DS). However, generalized version is applicable to largely arbitrary nonlinear PDE.

In the next section (Sec.\ref{Sec4}) we give general algorithm
deriving nonlinear $N$-wave type PDE. We introduce an 
 extra constrain allowing to write single  nonlinear PDE for single field. Characterization of solution space for derived nonlinear PDE is
given in Sec.\ref{IDE}. Sec.\ref{DR_DS} considers similar generalization of dressing method for DS.
Finally we represent some conclusions in Sec.\ref{conclusions}.

 \section{Derivation of multidimensional nonlinear $N$-wave equation}
\label{Sec4}
\label{DR_nl}

We start with usual  integral equation
\begin{eqnarray}
\label{Sec4:U}
\Phi(\lambda) =  \int\Psi(\lambda,\nu;x) U(\nu;x)
d\nu = \Psi(\lambda,\nu;x)*U(\nu;x)=\Psi*U,
\end{eqnarray}
where $*$ means integration over spectral parameter appearing
in both functions. 
There are two types of parameters in this equation. First, already mentioned {\it spectral parameters}
denoted by Greek letters $\lambda$, $\mu$, $\nu$ (for instance 
$\lambda=(\lambda_1,\dots,\lambda_{\dim \lambda})$), and, second, {\it additional 
parameters} denoted by $x$, 
$x= (x_1,\dots,x_{\dim x})$.  These additional parameters are 
{\it independent variables} of resulting nonlinear PDE. Besides, we reserve  
$k$   for scalar Fourier type parameter appearing in
integral representations of some functions. 
All  functions  are $Q\times Q$ matrices. 
We always assume $\dim x = M$ and $ \dim \lambda=  \dim \mu= \dim \nu=M+1$, where $M$ is dimensionality of resulting nonlinear PDE.
 
Eq.(\ref{Sec4:U}) is a linear equation for  
the {\it spectral function} $U(\lambda;x)$, where operator  
$\Psi(\lambda,\mu;x)*$ is required to be uniquely invertible, 
$\Phi$ is a diagonal matrix function
specified below. 
Integration is over whole space of vector spectral parameter $\nu$.
Function $\Psi(\lambda,\mu;x)$ is defined by the  following
formulae
 introducing $x$-dependence:
\begin{eqnarray}\label{Sec6:Fx}
\partial_{x_n}\Psi_{\alpha\beta}(\lambda,\mu;x)  +
\Big(h^n_\alpha(\lambda) +g^n_{\alpha\beta}(\mu)\Big)
\Psi_{\alpha\beta}(\lambda,\mu;x)  &=&
\Phi_\alpha(\lambda;x) B^{n}_\alpha
C_{\alpha\beta}(\mu;x),\\\nonumber
&&1\le n \le M.
\end{eqnarray}
 Here $C(\mu;x)$ is a new function, which will be characterized
 below; $h^n(\lambda)$ are diagonal
 and  $g^n(\mu)$ are arbitrary matrix
functions of argument;  $B^n$ are diagonal
constant matrices.
Short form of eq.(\ref{Sec6:Fx}) reads:
\begin{eqnarray}\label{Sec6:Fx2}
&& L^{n}_{\alpha\beta}(\lambda,\mu) \Psi_{\alpha\beta}(\lambda,\mu;x)= 
\Phi_\alpha(\lambda;x) B^{n}_\alpha
C_{\alpha\beta}(\mu;x),\;\;\\\nonumber
&& L^{n}_{\alpha\beta}(\lambda,\mu) (*) =
\Big(\partial_{x_{n}} +
 h^n_\alpha(\lambda)+  g^n_{\alpha\beta}(\mu) \Big)(*),\;\;1\le n
 \le M.
\end{eqnarray}
Remark that derivatives $\partial_{x_i}\Psi(\lambda,\mu;x)$ 
are not separated functions
of spectral parameters, unlike the $S$-integrable case
\cite{ZS,Z2}. 
Been overdetermined system of PDEs for function
$\Psi(\lambda,\mu;x)$, eqs. (\ref{Sec6:Fx2}) imply
compatibility conditions, which are following: 
\begin{eqnarray}\label{Sec4:comp1}
L^{n}_{\alpha\beta}(\lambda,\mu) \Big(\Phi_\alpha(\lambda;x)
B^{j}_\alpha C_{\alpha\beta}(\mu;x)\Big)&=&
L^{j}_{\alpha\beta}(\lambda,\mu) \Big(\Phi_\alpha(\lambda;x) B^{n}_\alpha
C_{\alpha\beta}(\mu;x)\Big),\;\;n\neq j.
\end{eqnarray}
  Without loss of generality,
 we put $j=1$, 
 and 
 \begin{eqnarray}
 \label{Sec1:simplification}
 B^1=I,\;\; h^1(\lambda)=g^1(\lambda)=0,
 \end{eqnarray}
 where $I$ is the unit matrix.
 Since each term in expanded form of eqs.(\ref{Sec4:comp1}) is  separated
 function of
parameters $\lambda$ and $\mu$, these equations  with $j=1$ 
are equivalent to two sets of 
 equations:
\begin{eqnarray}\label{Sec6:Phi_x}\label{Phi_x}
&& 
\partial_{x_n}\Phi(\lambda;x)  + h^n(\lambda) \Phi(\lambda;x)  - 
\partial_{x_1}\Phi (\lambda;x)
B^{n}=0,\;\;1<n\le M,
\\\label{c_x}
&& \partial_{x_n}C(\mu;x) +  
C^{1n}(\mu;x)  - 
B^{n}\partial_{x_1}C(\mu;x)  =
0,\;\;1<n\le M,\;\;
\end{eqnarray}
where
\begin{eqnarray}
\label{c^n}
C^{1n}_{\alpha\beta}(\mu;x) = C_{\alpha\beta}(\mu;x)
g^n_{\alpha\beta}(\mu).
\end{eqnarray}
Eqs.(\ref{Sec6:Phi_x}-\ref{c^n})  define
$\Phi$ and  $C$. We refer to functions $\Phi$, $C$ and $\Psi$ as {\it dressing functions}, where $\Psi$ is expressed in terms of $\Phi$ and $C$ due to eqs.(\ref{Sec6:Fx}).

Thus we have specified all functions appearing in 
eqs. (\ref{Sec4:U}) and (\ref{Sec6:Fx}). Now we demonstrate how
linear integral equation (\ref{Sec4:U}) is related with
appropriate multidimensional nonlinear PDE written for
fields expressible in terms of spectral function $U(\lambda;x)$ and
dressing functions.

System of nonlinear equations is generated by eq.(\ref{Sec6:Phi_x}).
  Derivation is very similar to derivation of  classical
  integrable equations \cite{ZS,Z2}. First of all, we use
 representation for $\Phi$ as  $\Psi*U$, see eq.(\ref{Sec4:U}). 
 Then using equations 
 (\ref{Sec6:Fx}) 
 for derivatives $\Psi_{x_n}$  
 we end up with homogeneous equations 
 in the form
 \begin{eqnarray}\label{Sec1:PsiEn}
 &&\Psi(\lambda,\mu;x)* E_n(\mu;x) =0,\;\;\\\nonumber
  &&  E_n(\lambda,x)= 
 U_{x_{n}}(\lambda,x) - U_{x_{1}}(\lambda,x) B^n  +
 U(\lambda,x) [B^{n},u(x)]
 - {\cal{G}}^n(\lambda,x),\;\;1<n\le M,
 \end{eqnarray} 
 where function $u$ is related with spectral function by
 the formula
 \begin{eqnarray}\label{Sec6:u}
 u(x)= C(\lambda,x)* U(\lambda,x)
 \end{eqnarray} 
and functions ${\cal{G}}^n$ satisfy the following equations:
\begin{eqnarray}\label{G}
\Psi(\lambda,\mu;x)* {\cal{G}}^n(\mu;x) &=& 
\Psi^n(\lambda,\mu;x)*
U(\mu;x),\;\;1<n\le M,\\\label{Psi^n}
\Psi^n_{\alpha\beta}(\lambda,\mu;x)&=&
\Psi_{\alpha\beta}(\lambda,\mu;x) \; g^n_{\alpha\beta}(\mu).
\end{eqnarray}
Later, function $u$ will be field of
nonlinear PDE.

Eqs.(\ref{G},\ref{Psi^n}) along with eq.(\ref{Sec4:U}) will be used  in
Sec.\ref{IDE} to analyse solution  space of nonlinear system. 
Inverting operator $\Psi*$ in 
eqs.(\ref{Sec1:PsiEn}) one gets
\begin{eqnarray}\label{Sec6:U2}
 E_n(\lambda;x):= U_{x_{n}}(\lambda;x) - U_{x_{1}}(\lambda;x) B^n  +
 U(\lambda;x) [B^{n},u(x)]   
 -{\cal{G}}^n(\lambda;x)  =0, \;\;1<n \le M.
 \end{eqnarray}
 In the case of  classical dressing method,  
 nonlinear integrable PDE can be received for function $u$
  applying $C(\lambda;x)*$ to (\ref{Sec6:U2}) 
and using eq.(\ref{c_x}) for
 $C_{x_n}$, $n>1$. Doing the same one gets in our case:
\begin{eqnarray}\label{eq:u}
 E^u_n(x)&:=& u_{x_{n}}(x) - u_{x_{1}}(x) B^n +u(x) [B^{n},u(x)]
 =\tilde
 {\cal{H}}^n(x) ,\\\nonumber
&& \tilde
 {\cal{H}}^n(x)=[B^n,u^1(x)] -C^{1n}(\mu;x)
 *U(\mu;x)    
 +C(\mu;x)*{\cal{G}}^n(\mu;x), \;\;1<n\le M,
\end{eqnarray}
where function $u^1$ is related with spectral function by the
formula similar to eq.(\ref{Sec6:u}):
\begin{eqnarray}
u^1(x)=C_{x_1}(\lambda;x)*U(\lambda;x).
\end{eqnarray}
Functions $u^1$ and  $\tilde{\cal{H}}^n$ are "intermediate"
functions which will be eliminated from the final system of
nonlinear PDEs.

System (\ref{eq:u}) has an obvious limit to classical
$(2+1)$-dimensional $S$-integrable  $N$-wave equation. In fact,
if $g^n=0$ for all $n$, (i.e. ${\cal{G}}^n(\lambda;x)=0$,
$\tilde{\cal{H}}^n(x)=[B^n,u^1(x)]$),
then  we may eliminate $u_1$
using two equations (\ref{eq:u}): $E^u_n$ and  $E^u_m$, $n\neq m$:
 \begin{eqnarray}\label{Sec0:Nw1}
[u_{x_n},B^m] - [u_{x_m},B^n] + B^m u_{x_1} B^n -  B^n u_{x_1}
B^m 
-[[u,B^m],[u,B^n]]=0.
\end{eqnarray}
This is the classical (2+1)-dimensional completely integrable $N$-wave equation, which has acceptable reduction $u_{\beta\alpha}=\bar u_{\alpha\beta}$, where bar means complex conjugation,
see, for instance, \cite{AC}.  
System 
(\ref{Sec6:U2}) with ${\cal{G}}^n=0$ becomes linear overdetermined system for eq.(\ref{Sec0:Nw1}),
 where  $U(\lambda;x)$ is a spectral function, i.e. eq.(\ref{Sec0:Nw1}) 
is compatibility condition for $E_n$ and $E_m$. This is well-known  common
feature of $S$-integrable models: they may be derived both
algebraically through compatibility condition of overdetermined
linear system and using dressing method.

 However, if $g^n \neq 0$ for all $n$, then  ${\cal{G}}^n\neq 0$ and system (\ref{Sec6:U2}) 
 may not be considered as a linear
 overdetermined system, since it has set of  spectral functions,
 such as $U(\lambda;x)$ and ${\cal{G}}^n(\lambda;x)$.
As a consequence, nonlinear eqs. (\ref{eq:u}) have  extra
functions
$\tilde{\cal{H}}^n(x)$ and may not be received 
as compatibility condition of the system (\ref{Sec6:U2})
through commutation of linear operators appearing in (\ref{Sec6:U2}).
So, similar to \cite{ZS}, the only way to derive system (\ref{eq:u}) from
eq.(\ref{Sec6:U2}) is the dressing
method.

The derived system (\ref{eq:u}) consists of ($M-1$)
equations and $M$ fields, which are $u$ and $\tilde {\cal{H}}^n$,
$1<n\le M$. In other words, it is not complete. 
In order to write a single nonlinear PDE for field $u$
we involve  another important deviation from the classical
approach. 

Let us  split $C(\lambda;x)$ into two factors:
\begin{eqnarray}\label{CGG}
&&C_{\alpha\beta}(\mu;x) =G^1_\alpha(\mu_1;x)
G^2_{\alpha\beta}(\mu;x), \\\nonumber
&& \partial_{x_n}G^1(\mu_1;x)  - 
B^{n}\partial_{x_1}G^1(\mu_1;x)  =
0,\;\;1<n\le M,\;\; \\\nonumber
&& \partial_{x_n}G^2(\mu;x) +  
G^{1n}(\mu;x)   =
0,\;\;1<n\le M,\;\;G^{1n}_{\alpha\beta}(\mu;x) = 
G_{\alpha\beta}(\mu;x) g^n_{\alpha\beta}(\mu), \\\nonumber
&& \partial_{x_1}G^2(\mu;x) =0
\end{eqnarray}
where eqs.(\ref{CGG}b-d) appear due to the eq.(\ref{c_x}).
Multiply eq.(\ref{Sec6:U2}) by $G^2(\lambda;x)$ from the left and integrate
over 
$\tilde \lambda= (\lambda_2,\dots,\lambda_{M+1})$.
One gets
\begin{eqnarray}\label{U2}
\tilde E_n(\lambda_1;x)&:=& \hat U_{x_{n}}(\lambda_1;x) - \hat
U_{x_{1}}(\lambda_1;x) B^n  +
 \hat U(\lambda_1;x) [B^{n},u(x)]   
 -\hat {\cal{F}}^n(\lambda_1;x)  =0, \;\;\\\nonumber
 && 1<n\le M,
 \end{eqnarray}
 where
 \begin{eqnarray}\label{hatU:def}
 \hat U(\lambda_1;x)&=& \int G^2(\lambda;x) U(\lambda;x) d\tilde
 \lambda,\\\nonumber
 \hat U^{1n}(\lambda_1;x)&=& 
 -\int G^2_{x_n}(\lambda;x) U(\lambda;x) d\tilde\lambda=
 \int G^{1n}(\lambda;x) U(\lambda;x) d\tilde\lambda
 ,
 \\\nonumber
 \hat {\cal{G}}^n(\lambda_1;x)&=& \int G^2(\lambda;x)
 {\cal{G}}^n(\lambda;x)d\tilde \lambda,\;\;
 \hat {\cal{F}}^n(\lambda_1;x) =  \hat {\cal{G}}^n(\lambda_1;x)-
 \hat U^{1n}(\lambda_1;x),\\\nonumber
 &&
 1<n\le M.
 \end{eqnarray}
We  will see in the next section that off-diagonal parts of 
$\hat U(\lambda_1;x) $ and
$\hat {\cal{G}}^n(\lambda_1;x) 
$ have arbitrary dependence on $x$. Thus we are able to introduce
one more relation among them. For instance, 
let 
\begin{eqnarray}\label{condition}
\sum_{i=2}^{M}S^i_{\alpha\beta}\Big(\hat 
{\cal{F}}^i_{\alpha\beta}(\lambda_1;x)
-\lambda_1\hat U_{\alpha\beta}(\lambda_1;x)
(B^i_\alpha-B^i_\beta)
\Big)
=0,\;\;\alpha\neq\beta,
\end{eqnarray}
where $S^i_{\alpha\beta}$ are constants.
Then eq.(\ref{U2}) gives ($\alpha\neq \beta$):
\begin{eqnarray}\label{lin:U}
\sum_{i=2}^{M}S^i_{\alpha\beta} 
\Big( \partial_{x_i} \hat U_{\alpha\beta}
(\lambda_1;x) - \partial_{x_1}\hat U_{\alpha\beta}(\lambda_1;x)
B^i_\beta  +
 \sum_{\gamma=1}^Q\hat U_{\alpha\gamma}(\lambda_1;x)
 u_{\gamma\beta}(x) (B^{i}_\gamma-B^{i}_\beta) -&&
 \\\nonumber
 \lambda_1\hat U_{\alpha\beta}(\lambda_1;x)
(B^i_\alpha-B^i_\beta)\Big) &=&0, \\\label{lin:Ub}
 \sum_{i=2}^{M}
S^i_{\alpha\beta}(B^i_\alpha-B^i_\beta) &=&0.
\end{eqnarray}
This equation is a linear equation for the spectral function $\hat U^{of}$; additional relation (\ref{lin:Ub}) is introduced to eliminate diagonal part of $\hat U$ from the nonlinear term of (\ref{lin:U}). 

Multiply  this equation by $G^1_\alpha(\lambda_1;x)$ from the left, integrate over
$\lambda_1$ and assume that $G^1_{x_1}(\lambda_1;x)=\lambda_1 G^1(\lambda_1;x)$:
\begin{eqnarray}\label{nl}
\sum_{i=2}^{M}S^i_{\alpha\beta} 
\Big( \partial_{x_i} u_{\alpha\beta}
(x) - \partial_{x_1}u_{\alpha\beta}(x)
B^i_\beta  +
 \sum_{{\gamma=1}\atop{\gamma\neq\alpha\neq\beta}}^Q u_{\alpha\gamma}(x)
 u_{\gamma\beta}(x) (B^{i}_\gamma-B^{i}_\beta) \Big) =0,\;\;\alpha\neq\beta
\end{eqnarray}
which becomes  $N$-wave equation if, along with (\ref{lin:Ub}), one requires
\begin{eqnarray}
S^i_{\alpha\beta}  =
  S^i_{\beta\alpha},\;\;
  u_{\beta\alpha}=\bar u_{\alpha\beta}.
\end{eqnarray}
Thus, nonlinear eq.(\ref{nl}) is equivalent to linear eq.(\ref{lin:U}) where spectral function $\hat U^{of}(\lambda_1;x)$ is related with dressing functions by the eqs.(\ref{Sec4:U}-\ref{c^n},\ref{CGG},\ref{hatU:def}). Detailed discussion of this relation is represented in the next subsection.

\subsection{Analysis of the system (\ref{Sec4:U}-\ref{c^n},\ref{CGG},\ref{hatU:def},\ref{lin:U})}
\label{IDE}
\label{Sec63}

In this section we characterize solution space of nonlinear equation
(\ref{nl}) in terms of dressing functions $\Psi$, $\Phi$ and $C$.  
First step is solving equations
(\ref{Sec6:Fx},\ref{Sec6:Phi_x},\ref{c_x}) 
for $\Psi(\lambda,\mu;x)$, $\Phi(\lambda;x)$ and $C(\mu;x)$. 
Eq.(\ref{Sec6:Fx}) is nonhomogeneous equation 
for $\Psi(\lambda,\mu;x)$, so we take the following solution:
\begin{eqnarray}\label{Sec5:F}
\Psi_{\alpha\beta}(\lambda,\mu;x) &=& \partial_{x_1}^{-1}
\Big(\Phi_\alpha(\lambda;x) 
C_{\alpha\beta}(\mu;x)\Big) +
\delta_{\alpha\beta}\delta(\lambda-\mu) 
e^{-\sum\limits_{j=2}^M \Big(h^j_\alpha(\lambda) +
g^j_{\alpha\beta}(\mu)\Big)},\\\nonumber
&&
\delta(\lambda-\mu)= \prod_{i=1}^{M+1}\delta(\lambda_i-\mu_i)
\end{eqnarray}
(remember that dimension of spectral parameters is $M+1$),
where $\delta_{\alpha\beta}$ is Kronecker delta symbol, first term is a particular solution of nonhomogeneous equation, while the second term is particular solution of homogeneous equation associated with eq.(\ref{Sec6:Fx}).
Function (\ref{Sec5:F}) is not general solution of (\ref{Sec6:Fx}), but this is enough for
our algorithm.

Solutions of eqs. (\ref{Sec6:Phi_x},\ref{c_x})
in view of (\ref{CGG}) read
\begin{eqnarray}\label{Sec4:chi_sol_t}
\label{Phi_def}
\Phi_{\alpha}(\lambda;x)&=& 
\int \Phi^0_{\alpha}(\lambda,k) 
e^{K^\Phi_{\alpha}(\lambda,k;x)} d k ,\;\;
K^\Phi_{\alpha}(\lambda,k;x)=k x_1+  \sum_{j=2}^{M}
\left(k B^{j}_\alpha  - h^j_\alpha(\lambda)  \right)
 x_j,\\\label{Sec4:c_sol}
 C_{\alpha\beta}(\mu;x)&=& G^1_\alpha(\mu_1;x) G^2_{\alpha\beta}(\mu;x)
 ,\\\nonumber
G^1_\alpha(\mu_1;x) &=&
e^{K^{G^1}_\alpha(\mu_1;x)},\;\;
G^2_{\alpha\beta}(\mu;x)= 
 e^{K^{G^2}_{\alpha\beta}(\mu;x)} C^0_{\alpha\beta}(\mu) 
,\\\nonumber
&&
K^{G^1}_\alpha(\mu_1;x)=\mu_1\left(x_1 + \sum_{i=2}^{M} B^i_\alpha
x_i \right) ,\;\;
K^{G^2}_{\alpha\beta}(\mu;x)=-\sum_{j=2}^{M}
     g^j_{\alpha\beta}(\mu)  x_{j},
 \end{eqnarray}
where  parameter $k$  is scalar.

Hereafter we take
\begin{eqnarray}\label{Sec4:cphi}
\Phi^0(\lambda,k) = 
\delta(\lambda_2-k) I.
\end{eqnarray}
 Thus expression
 (\ref{Sec5:F})  may be written in  explicit
 form:
\begin{eqnarray}\label{Sec6:Fexpl}
\Psi_{\alpha\beta}(\lambda,\mu;x) &=&  
\frac{e^{K^\Phi_{\alpha}(\lambda,\lambda_2;x)+
K^{G^1}_\alpha(\mu_1;x)+K^{G^2}_{\alpha\beta}(\mu;x)} C^0_{\alpha\beta}(\mu) 
 }{\lambda_2+\mu_1} 
 + 
\delta_{\alpha\beta}\delta(\lambda-\mu)
 e^{-\sum\limits_{j=2}^M \Big(h^j_\alpha(\lambda) +
g^j_{\alpha\beta}(\mu)\Big)}.
\end{eqnarray}
Due to the last term in eq.(\ref{Sec6:Fexpl}), eq.(\ref{Sec4:U}) has term $e^{-\sum\limits_{j=2}^M\Big( h^j_\alpha(\lambda)+ 
g^j_{\alpha\alpha}(\lambda)\Big) x_j} U_{\alpha\beta}(\lambda;x)$. However, we would like to eliminate factor ahead of $U$ in this term for convenience of subsequent constructions. 
To do this, we 
 multiply 
  eqs.(\ref{Sec4:U},\ref{G}) by $e^{\sum\limits_{j=2}^M\Big( h^j_\alpha(\lambda) + 
g^j_{\alpha\alpha}(\lambda)\Big) x_j}$:
\begin{eqnarray}\label{Sec4:E0^U}
 E^U(\lambda;x)
&:=& U(\lambda;x)= - 
\partial^{-1}_{x_1} \Big(\Phi^1(\lambda;x) C(\mu;x)\Big)*U(\mu;x) +
 \Phi^1(\lambda;x), \\
\label{G2_2}
E^{G^n}(\lambda;x) &:=& {\cal{G}}^n(\lambda;x) =  - 
\partial^{-1}_{x_1} \Big(\Phi^1(\lambda;x) 
C(\mu;x)\Big)*{\cal{G}}^n(\mu;x)
+
\\\nonumber
&&
 \partial^{-1}_{x_1} \Big(\Phi^1(\lambda;x) 
C^{1n}(\mu;x)\Big)*
U(\mu;x) +
  U^n(\lambda;x),\;\;n>1 ,
\end{eqnarray}
where
\begin{eqnarray}\label{Phi^1} 
\Phi^1_{\alpha}(\lambda;x)&=&e^{\sum\limits_{j=2}^{M}
\Big(h^j_\alpha(\lambda) + g^j_{\alpha\alpha}(\lambda)\Big)}
 \Phi_{\alpha}(\lambda;x)= 
e^{K^{\Phi^1}_{\alpha}(\lambda;x)}  ,\;\;\\\nonumber
&&
K^{\Phi^1}_{\alpha}(\lambda;x)= \lambda_2 x_1+  \sum_{j=2}^{M}
\left(\lambda_2 B^{j}_\alpha  + g^j_{\alpha\alpha}(\lambda)  \right)
 x_j,\\
 \label{U^n}
   U^n_{\alpha\beta} (\lambda;x)&=&
   g^n_{\alpha\alpha}(\lambda) U_{\alpha\beta}(\lambda;x).
\end{eqnarray}
Below we need function
\begin{eqnarray}
G^{1n}_{\alpha\beta}(\lambda;x) = G^2_{\alpha\beta}(\lambda;x) g^n_{\alpha\beta}(\lambda).
\end{eqnarray}
Applying $\int  d\tilde \lambda G^2(\lambda;x) \cdot$ to eqs 
(\ref{Sec4:E0^U},\ref{G2_2}) and $\int  d\tilde \lambda G^{1n}(\lambda;x) \cdot$
to (\ref{Sec4:E0^U})  one gets equations for $\hat U$, $\hat {\cal{G}}^n$ and
$\hat U^n$:
\begin{eqnarray}\label{hatU}
 \hat U(\lambda_1;x)&=& - 
\int\partial^{-1}_{x_1} \Big(\hat\Phi(\lambda_1;x)
G^1(\mu_1;x)\Big)\hat U(\mu_1;x)d\mu_1 +
 \hat\Phi(\lambda_1;x), 
\\\label{hatG2_2}
 \hat{\cal{G}}^n(\lambda_1;x) &=&  - 
\int\partial^{-1}_{x_1} \Big(\hat\Phi(\lambda_1;x) 
G^1(\mu_1;x)\Big)\Big(\hat{\cal{G}}^n(\mu_1;x) 
-
\hat U^{1n}(\mu_1;x)\Big) d\mu_1 +
\\\nonumber
&&
  \hat U^{2n}(\lambda_1;x),\;\;1<n\le M , 
  \\
\label{hatG^n}
  \hat U^{1n}(\lambda_1;x)&=& - 
\int\partial^{-1}_{x_1} \Big(\hat\Phi^{1n}(\lambda_1;x)
G^1(\mu_1;x)\Big)\hat U(\mu_1;x) d\mu_1 +
 \hat\Phi^{1n}(\lambda_1;x),\;\; 1<n\le M,
\end{eqnarray}
where
\begin{eqnarray}
&&
\hat \Phi(\lambda_1;x)=\int G^2(\lambda;x) \Phi^1(\lambda;x) 
d\tilde\lambda,\;\;\hat \Phi^{1n}(\lambda_1;x)=\int G^{1n}(\lambda;x) \Phi^1(\lambda;x) 
d\tilde\lambda,\\\nonumber
&&
\hat U^{2n}(\lambda_1;x)=
\int G^2(\lambda;x) U^{n}(\lambda;x) d\tilde\lambda=
\int G^{2n}(\lambda;x) U(\lambda;x) d\tilde\lambda
,
\\\nonumber
&&
G^{2n}_{\alpha\beta}(\lambda;x)= G^2_{\alpha\beta}(\lambda;x) g^n_{\beta\beta}(\lambda).
\end{eqnarray}
Equation for $\hat U^{2n}$  follows from eq.(\ref{Sec4:E0^U}) after applying 
$\int  d\tilde \lambda G^{2n}(\lambda;x) \cdot$:
\begin{eqnarray}\label{hatU^2n}
\hat U^{2n}(\lambda_1;x)&=& - 
\int\partial^{-1}_{x_1} \Big(\hat\Phi^{2n}(\lambda_1;x) G^1(\mu_1;x)\Big) \hat U(\mu_1;x)d\mu_1 +
 \hat\Phi^{2n}(\lambda_1;x),\;\;n>1,\\\nonumber
 &&
 \hat \Phi^{2n}(\lambda_1;x)=\int G^{2n}(\lambda;x) \Phi^1(\lambda;x) d\tilde\lambda.
\end{eqnarray}

By construction, 
function $\hat \Phi(\lambda_1;x)$ has arbitrary dependence on variables $x$, if, for instance, $g^i_{\alpha\beta}(\lambda)= \lambda_{i+1} \hat g^i_{\alpha\beta}$, where 
$\hat g^i_{\alpha\beta}$ are constants, $i=2,\dots,M$. Due to this fact $\hat \Phi(\lambda_1;x)$ may solve
  equation (\ref{condition}). Let us transform eq.(\ref{condition})
substituting eqs.(\ref{hatU}-\ref{hatU^2n}):
\begin{eqnarray}
\label{condition_f}
\sum_{i=2}^{M}S^i_{\alpha\beta}\Big\{
\partial_{x_i}\hat \Phi_{\alpha\beta}(\lambda_1;x)-\partial_{x_1} \hat \Phi_{\alpha\beta}(\lambda_1;x) B^i_\beta - \lambda_1 \hat \Phi_{\alpha\beta}(\lambda_1;x)(B^i_\alpha- B^i_\beta)-&&\\\nonumber
\int\sum_{\gamma=1}^Q\Big[
\partial_{x_1}^{-1} \Big(\hat 
\Phi_{\alpha\gamma}(\lambda_1;x) G^1_{\gamma}(\mu_1;x)\Big) \hat{\cal{F}}^i_{\gamma\beta}(\mu_1;x)+&&\\\nonumber
\partial_{x_1}^{-1} \Big(\big(\partial_{x_i}\hat \Phi_{\alpha\gamma}(\lambda_1;x)-\partial_{x_1} \hat \Phi_{\alpha\gamma}(\lambda_1;x) B^i_\gamma\big)G^1_{\gamma}(\mu_1;x)\Big) \hat U_{\gamma\beta}(\mu;x)-
&&\\\nonumber
\lambda_1 \partial_{x_1}^{-1} \Big(\hat \Phi_{\alpha\gamma}(\lambda_1;x) G^1_{\gamma}(\mu_1;x)\Big) \hat U_{\gamma\beta}(\mu_1;x)
(B^i_\alpha- B^i_\beta) 
\Big]d\mu_1 \Big\} &=& 0
\end{eqnarray}
where $\alpha\neq\beta$  
and eqs.(\ref{hatU}-\ref{hatU^2n}) give us
\begin{eqnarray}\label{U^23n}
\hat{\cal{F}}^n(\lambda_1;x) &=&  - 
\int\partial^{-1}_{x_1} \Big(\hat\Phi(\lambda_1;x) 
G^1(\mu_1;x)\Big)\hat{\cal{F}}^n(\mu_1;x) d\mu_1
 -  \\\nonumber
 &&
\int\partial^{-1}_{x_1} \Big(\hat\Phi_{x_n}(\lambda_1;x)-\hat\Phi_{x_1}(\lambda_1;x)B^n
\Big) G^1(\mu_1;x) \hat U(\mu_1;x) d\mu_1 
 +  \\\nonumber
 &&
 \hat\Phi_{x_n}(\lambda_1;x)-\hat\Phi_{x_1}(\lambda_1;x)B^n,\;\;
 1<n\le M.
\end{eqnarray}
Deriving eqs.(\ref{condition_f},\ref{U^23n}), we took into account an obvious relation
\begin{eqnarray}
\hat\Phi^{2n}_{\alpha\beta} - \hat \Phi^{1n}_{\alpha\beta}  = 
\partial_{x_n} \hat \Phi_{\alpha\beta}-\partial_{x_1} 
\hat \Phi_{\alpha\beta} B^n_\beta, \;\;\alpha\neq \beta, \;\; 1<n\le M .
\end{eqnarray}
Note, that diagonal elements of $\hat\Phi$,
\begin{eqnarray}
\label{diag}
\hat \Phi_{\alpha\alpha}(\lambda_1;x) =\int C^0_{\alpha\alpha}(\lambda) e^{\lambda_2 \Big(x_1 + \sum_{\i=2}^M B^i_\alpha x_i\Big)}d\tilde \lambda,
\end{eqnarray}
may be arbitrary functions of single independent variable and
$\hat\Phi^{2n}_{\alpha\alpha} - \hat \Phi^{1n}_{\alpha\alpha} =0$.

System (\ref{hatU},\ref{condition_f},\ref{U^23n})
represent a complete nonlinear system of equations allowing to find $\hat U(\lambda_1;x)$ and $\hat\Phi(\lambda_1;x)$. Since $u(x) = \int G^1(\lambda_1;x) \hat U(\lambda_1;x) d\lambda_1$, this system is alternative form of the nonlinear equation (\ref{nl}). 
In particular case $C^0(\mu)=\delta(\mu_1) \tilde C^0(\tilde \mu)$, one has $G^1(0;x)=I$, and this system reduces to PDE for $\varphi(\lambda_1;x)= \partial_{x_1}^{-1} \hat \Phi(\lambda_1;x)$ (below $\alpha\neq\beta$):
\begin{eqnarray}\label{RhatU}
 \hat U(\lambda_1;x)= - 
\varphi(\lambda_1;x)
u(x) +
 \partial_{x_1}\varphi(\lambda_1;x),\;\;\; u(x)=\hat U(0;x),\hspace{3cm}&& \\
\label{condition_f2}
\sum_{i=2}^{M}S^i_{\alpha\beta}\Big\{\partial_{x_1}\Big(
\partial_{x_i}\varphi_{\alpha\beta}(\lambda_1;x)-\partial_{x_1} \varphi_{\alpha\beta}(\lambda_1;x) B^i_\beta-\lambda_1\varphi_{\alpha\beta}(\lambda_1;x)(B^i_\alpha- B^i_\beta)\Big)-&&\\\nonumber
\sum_{\gamma=1}^Q\Big[
\varphi_{\alpha\gamma}(\lambda_1;x)  \hat{\cal{F}}^i_{\gamma\beta}(0;x)+&&\\\nonumber \Big(\partial_{x_i}\varphi_{\alpha\gamma}(\lambda_1;x)-\partial_{x_1} \varphi_{\alpha\gamma}(\lambda_1;x) B^i_\gamma\Big) u_{\gamma\beta}(x)-
\lambda_1 \varphi_{\alpha\gamma}(\lambda_1;x) u_{\gamma\beta}(x)
(B^i_\alpha- B^i_\beta) 
\Big] \Big\} &=& 0,
\end{eqnarray}
\begin{eqnarray}
\label{RU^23n}
 \hat{\cal{F}}^n(\lambda_1;x)&=&  - 
\varphi(\lambda_1;x) 
\hat{\cal{F}}^n(0;x) 
-\\\nonumber
&&
\Big(\varphi_{x_n}(\lambda_1;x)-
\varphi_{x_1}(\lambda_1;x)B^n\Big)   u(x) 
+ \varphi_{x_1x_n}(\lambda_1;x)-
\varphi_{x_1x_1}(\lambda_1;x)B^n
,\\\nonumber
&&
1<n\le M,
\end{eqnarray}

Eqs. (\ref{RhatU},\ref{RU^23n}) with $\lambda_1=0$ give us
\begin{eqnarray}\label{potentials}
u(x)&=&\Big(I+\varphi(0;x)\Big)^{-1}\varphi_{x_1}(0;x),
\\\nonumber
\hat{\cal{F}}^n(0;x) 
 &=& 
\Big(I+\varphi(0;x)\Big)^{-1}
\Big[\varphi_{x_1x_n}(0;x)-
\varphi_{x_1x_1}(0;x)B^n-\\\nonumber
&&
\Big(\varphi_{x_n}(0;x)-
\varphi_{x_1}(0;x)B^n\Big)   u(x) 
\Big]
\end{eqnarray}
We see that eq.(\ref{condition_f2}) is linear PDE for $\varphi^{of}(\lambda_1;x)$, $\lambda_1\neq 0$,  with "boundary" function $\varphi(0;x)$ satisfying (\ref{potentials}a).
By construction, if $\lambda_1=0$, then  eq.(\ref{condition_f2}) is projected  into (\ref{nl}), i.e. calculation of evolution of  $\varphi^{of}(0;x)$
is equivalent to solving original nonlinear PDE (\ref{nl}). However, from another point of view, this evolution may be found as  $\lim\limits_{\lambda_1\to 0} \varphi^{of}(\lambda_1;x)$.

The simplest algorithm for numerical construction of particular solutions to (\ref{nl}) is following.

For given arbitrary $\varphi(\lambda_1;x)|_{x_M=0}$ we find  
$u(x)|_{x_M=0}$ and $\hat{\cal{G}}^n(0;x) |_{x_M=0}
- \hat U^{1n}(0;x)|_{x_M=0}$  using   (\ref{potentials}). Then solve (\ref{condition_f2}) for $\varphi^{of}_{x_M}(\lambda_1;x)|_{x_M=0}$. 
Using Tailor formulae we approximate $\varphi^{of}(\lambda_1;x)|_{x_M=\Delta t}$ :
\begin{eqnarray}
\varphi^{of}(\lambda_1;x)|_{x_M=\Delta t} \approx 
\varphi^{of}(\lambda_1;x)|_{x_M=0} + \Delta t \; \varphi^{of}_{x_M}(\lambda_1;x)|_{x_M=0}.
\end{eqnarray}
Evolution of diagonal elements $\varphi_{\alpha\alpha}(\lambda_1;x)$ is fixed by $\varphi_{\alpha\alpha}(\lambda_1;x)|_{x_M=0}$  due to (\ref{diag}).
Substitute this result into (\ref{potentials}) we find
$u(x)|_{x_M=\Delta t}$ and $\hat{\cal{G}}^n(0;x) |_{x_M=\Delta t}
- \hat U^{1n}(0;x)|_{x_M=\Delta t}$.
Then eq.(\ref{condition_f2}) gives 
$\varphi^{of}_{x_M}(\lambda_1;x)|_{x_M=\Delta t}$, and so on. 
Solving the Initial Value Problem (IVP) (i.e. construction of $u^{of}(x)$ for given initial data $u^{of}(x)|_{x_M=0}$) is more complicated and will not be considered here, since it seems to be not simpler then direct numerical solving of IVP for (\ref{nl}). 

Let us remark in the end of this section, that eq.(\ref{condition}) is not the only admittable constrain. Instead of zero in the rhs of this equation one might use expression  $L\big(\hat U^{of}(\lambda_1;x\big); u^{of}(x))$ which is linear differential operator applied to $\hat U^{of}(\lambda_1;x)$.  Coefficients of this operator  depend on field $u^{of}(x)$ and its derivatives. Then expression $L\big(\hat U^{of}(\lambda_1;x\big); u^{of}(x))$  appears in the rhs of (\ref{lin:U}). The only requirement to $L$ is that 
after multiplying eq.(\ref{lin:U}) by  $G^1_\alpha(\lambda_1;x)$ 
from the left and integrating over $\lambda_1$ one gets nonlinear PDE for  $u^{of}$. This new PDE (which replaces eq.(\ref{nl})) may be  
largely arbitrary nonlinear PDE for $u^{of}$. So, as for now, represented  multidimensional version of the dressing method is not the method for solving of nonlinear PDE, but it gives a new  representation of nonlinear PDE. This situation is equivalent to the situation appearing when Fourier method is applied to PDE other then linear PDE with constant coefficients.

\section{Derivation of multidimensional Nonlinear Shr\"odinger Equation}
\label{DR_DS}

In the previous section we demonstrated that (largely) arbitrary nonlinear PDE can be transformed using a variant of multidimensional generalization of the dressing method for (2+1)-dimensional $N$-wave equation.  
In this section we show that similar construction may be performed  starting with the dressing method for (2+1)-dimensional DS. We use notations of the Sec.\ref{DR_nl}. For simplicity,  we  take $Q=2$, i.e. consider $2\times 2$ matrix equations. 
 Variables $x_i$ are introduced by the following system:
\begin{eqnarray}\label{DS_Sec6:Fx}
\partial_{x_n}\Psi_{\alpha\beta}(\lambda,\mu;x)  +
\Big(h^n_\alpha(\lambda) +g^n_{\alpha\beta}(\mu)\Big)
\Psi_{\alpha\beta}(\lambda,\mu;x)  &=&
\Phi_\alpha(\lambda;x) B^{n}_\alpha
C_{\alpha\beta}(\mu;x),\\\nonumber
&& 1\le n < M\\\nonumber
\partial_{x_M}\Psi_{\alpha\beta}(\lambda,\mu;x)  +
\Big(h^M_\alpha(\lambda) +g^M_{\alpha\beta}(\mu)\Big)
\Psi_{\alpha\beta}(\lambda,\mu;x)  &=&
\partial_{x_1}\Phi_\alpha(\lambda;x) B^{M}_\alpha
C_{\alpha\beta}(\mu;x)-\\\nonumber
&&
\Phi_\alpha(\lambda;x) B^{M}_\alpha
\partial_{x_1}C_{\alpha\beta}(\mu;x),
\end{eqnarray}
where the first equation is identical to (\ref{Sec6:Fx}).
Since $Q=2$, only two $B^i$ are linearly independent, so we may put $B^i=0$, $i>2$ without loss of generality. Let, in addition, $B^1=I$, $B^M=B^2={\mbox{diag}}(1,-1)$, $h^1=g^1=0$. 
Compatibility of (\ref{DS_Sec6:Fx})  results in (compare with Sec.(\ref{DR_nl})):
\begin{eqnarray}\label{DS_Sec6:Phi_x}\label{DS_Phi_x}
&&
\partial_{x_{2}}\Phi(\lambda;x)  + h^2(\lambda) \Phi(\lambda;x)  - 
\partial_{x_1}\Phi (\lambda;x) B^2=0,\\\nonumber
&&
 \partial_{x_{n}}\Phi(\lambda;x)  + h^{n}(\lambda) \Phi(\lambda;x)  =0,\;\;2<n<M
\\\nonumber
&&
\partial_{x_{M}}\Phi(\lambda;x)  + h^{M}(\lambda) \Phi(\lambda;x)  - 
\partial^2_{x_1}\Phi (\lambda;x)B^2
=0,
\\\nonumber\\
\label{DS_c_x}
&&\partial_{x_{2}}C(\mu;x) +  
C^{12}(\mu;x)  -
B^2\partial_{x_1}C(\mu;x)  =
0,\\\nonumber
&&
\partial_{x_{n}}C(\mu;x) +  
C^{1n}(\mu;x)    =
0,\;\;2<n<M,\\\nonumber
&&
\partial_{x_{M}}C(\mu;x) +  
C^{1M}(\mu;x)  + 
B^2\partial^2_{x_1}C(\mu;x)  =
0,
\end{eqnarray}
where
\begin{eqnarray}
\label{DS_c^n}
C^{1n}_{\alpha\beta}(\mu;x) = C_{\alpha\beta}(\mu;x)
g^{n}_{\alpha\beta}(\mu), \;\;1<n \le M.
\end{eqnarray}
Eqs.(\ref{DS_Sec6:Phi_x}-\ref{DS_c^n})  define
$\Phi$ and  $C$.

System of nonlinear equations is generated by eq.(\ref{DS_Sec6:Phi_x}).
  Derivation is very similar to derivation carried out in  Sec.\ref{DR_nl}. First of all, we use
 representation for $\Phi$ as  $\Psi*U$, see eq.(\ref{Sec4:U}). 
 Then using equations 
 (\ref{DS_Sec6:Fx}) 
 for derivatives $\Psi_{x_n}$ and inverting $\Psi*$   
 we end up with system of linear equations
 in the form
 \begin{eqnarray}\label{DS_U2_1}
  E_2(\lambda;x)&:=& U_{x_{2}}(\lambda;x) - U_{x_{1}}(\lambda;x) B^2  +
 U(\lambda;x) [B^{2},u(x)]   
 -{\cal{G}}^2(\lambda;x)  =0, \\\label{DS_U2_2}
  E_n(\lambda;x)&:=& U_{x_{n}}(\lambda;x)    
 -{\cal{G}}^n(\lambda;x)  =0, \\\nonumber
 &&
 2<n < M\\\label{DS_U2_3}
    E_{M}(\lambda,x)&:=& 
 U_{x_{M}}(\lambda,x) - U_{x_1 x_{1}}(\lambda,x) B^2  +
  U(\lambda,x) \Big(u(x)[B^2,u(x)] -
  \\\nonumber
  && 2 u_{x_1}(x) B^2 +[u^1,B^2]\Big) + U_{x_1}(\lambda,x) [B^2,u] - {\cal{G}}^{M}(\lambda;x)=0
 \end{eqnarray} 
 where functions $u$ and $u^1$ are related with spectral functions by
 the formula
 \begin{eqnarray}\label{DS_Sec6:u}
 u(x)= C(\lambda,x)* U(\lambda,x),\;\; u^1(x)= C_{x_1}(\lambda,x)* U(\lambda,x)
 \end{eqnarray} 
and functions ${\cal{G}}^n$ satisfy the following equations:
\begin{eqnarray}\label{DS_G}
\Psi(\lambda,\mu;x)* {\cal{G}}^n(\mu;x) &=& 
\Psi^n(\lambda,\mu;x)*
U(\mu;x),\\\label{DS_Psi^n}\nonumber
\Psi^n_{\alpha\beta}(\lambda,\mu;x)&=&
\Psi_{\alpha\beta}(\lambda,\mu;x) \; g^n_{\alpha\beta}(\mu),\;\;1<n\le M.
\end{eqnarray}
Later, function $u$ will be  field in the 
nonlinear PDE.
 
 In the case of  classical dressing method,  
 nonlinear integrable PDE can be received for function $u$
  applying $C(\lambda;x)*$ and $C_{x_1}(\lambda;x)*$ to (\ref{DS_U2_1}),  applying $C(\lambda;x)*$ to (\ref{DS_U2_3}) 
and using eqs.(\ref{DS_c_x}) for
 $C_{x_n}$, $n>1$. Doing the same one gets in our case:
\begin{eqnarray}\label{DS_eq:u}
 E^u_{02}(x)&:=& u_{x_{2}}(x) - u_{x_{1}}(x) B^2 +u(x) [B^{2},u(x)]
 =\tilde
 {\cal{H}}^{02}(x) ,\\\label{DS_eq:u1}
  E^u_{12}(x)&:=& u^1_{x_{2}}(x) - u^1_{x_{1}}(x) B^2 +u^1(x) [B^{2},u(x)]
 =\tilde
 {\cal{H}}^{12}(x) ,\\\label{DS_eq:uM}
E^u_M(x)&:=& u_{x_{M}}(x) - u_{x_1 x_{1}}(x) B^2 +u(x) \Big(u(x)[B^{2},u(x)]-2u_{x_1} B^2\Big) + \\\nonumber
&&u_{x_1}(x) [B^2,u(x)] 
 =\tilde
 {\cal{H}}^{0M}(x) ,
 \\\nonumber 
&& \tilde
 {\cal{H}}^{02}(x)=[B^2,u^1(x)] -C^{12}(\mu;x)
 *U(\mu;x)    
 +C(\mu;x)*{\cal{G}}^2(\mu;x),\\\nonumber 
&& \tilde
 {\cal{H}}^{12}(x)=[B^2,u^2(x)] -C^{12}_{x_1}(\mu;x)
 *U(\mu;x)    
 +C_{x_1}(\mu;x)*{\cal{G}}^2(\mu;x),\\\nonumber
&& \tilde
 {\cal{H}}^{0M}(x)=- 2u^1_{x_1}(x) B^2 -u(x)[u^1(x),B^2] + u^1(x)[B^2,u(x)]-[B^2,u^2(x)] -\\\nonumber
 &&
 C^{1M}(\mu;x)
 *U(\mu;x)    
 +C(\mu;x)*{\cal{G}}^M(\mu;x), 
\end{eqnarray}
where function $u^2$ is related with spectral function by the
formula similar to eq.(\ref{Sec6:u}):
\begin{eqnarray}
u^2(x)=C_{x_1x_1}(\lambda;x)*U(\lambda;x).
\end{eqnarray}
Remark that eqs.(\ref{DS_eq:u}) coincide with (\ref{eq:u}) 
where $n=2$ and $\tilde{\cal{H}}^2=\tilde{\cal{H}}^{02}$.
Functions $u^1$, $u^2$  and  $\tilde{\cal{H}}^{in}$ are "intermediate"
functions which will be eliminated from the final 
nonlinear PDE.

System (\ref{DS_eq:u}-\ref{DS_eq:uM}) has an obvious limit to classical
$(2+1)$-dimensional $S$-integrable  
DS. In fact,
if $g^2=g^M=0$, (i.e. ${\cal{G}}^2(\lambda;x)={\cal{G}}^M(\lambda;x)=0$,
$\tilde{\cal{H}}^{02}(x)=[B^2,u^1(x)]$,
$\tilde{\cal{H}}^{12}(x)=[B^2,u^2(x)]$, $\tilde{\cal{H}}^{0M}(x)=- 2u^1_{x_1}(x) B^2 -u(x)[u^1(x),B^2] + u^1(x)[B^2,u(x)]-[B^2,u^2(x)]$),
then  we may eliminate $u_1$ and $u_2$ from eq.(\ref{DS_eq:uM})
using  equations (\ref{DS_eq:u}) and  (\ref{DS_eq:u1}):
 \begin{eqnarray}\label{Sec0:DS}
{\cal{E}}&:=&[u^{of}_{x_{M+1}},\sigma]-u^{of}_{x_1x_1}
-u^{of}_{x_2x_2}
-8 u_{12} u_{21} u^{of} -4 \varphi u^{of}
=0\\\label{Sec1:varphi}
 &&\varphi_{x_2 x_2} - \varphi_{x_1x_1} = 4
 (u_{12} u_{21})_{x_1x_1},\;\;
\varphi=(u_{11}+u_{22})_{x_1},
\end{eqnarray}
where 
\begin{eqnarray}
u=\left(
\begin{array}{cc}
u_{11} & u_{12} \cr
u_{21} & u_{22}
\end{array}
\right),
\end{eqnarray}
which is DS if $x_{M}=i t$, $i^2=-1$, $u_{21}=\bar u_{12}$.
Eqs.
(\ref{DS_U2_1})  and (\ref{DS_U2_3}) 
with ${\cal{G}}^2={\cal{G}}^M=0$ become linear overdetermined system for
this equation where  spectral function is $U(\lambda;x)$, i.e. eq.(\ref{Sec0:DS}) 
is compatibility condition for $E_2$ and $E_M$. 

 However, if $g^2 \neq 0$ and $g^M \neq 0$, then system (\ref{DS_U2_1},\ref{DS_U2_3}) 
 may not be considered as a linear
 overdetermined system, since it has set of spectral functions,
 such as $U(\lambda;x)$ and ${\cal{G}}^n(\lambda;x)$.
As a consequence, nonlinear eqs. (\ref{DS_eq:u}-\ref{DS_eq:uM}) have  extra
fields
$\tilde{\cal{H}}^{in}(x)$, $i=0,1$, and may not be received 
as compatibility condition of the system (\ref{DS_U2_1},\ref{DS_U2_3})
through commutation of linear operators.
So, similar to Sec.\ref{DR_nl}, the only way to derive nonlinear system (\ref{DS_eq:u}-\ref{DS_eq:uM}) from
eqs.(\ref{DS_U2_1},\ref{DS_U2_3}) is the dressing
method.

Similar to Sec.\ref{IDE}, we can take largely arbitrary equation for $\hat U^{of}(\lambda_1;x)$ resulting to largely arbitrary nonlinear PDE for field $u^{of}$. 
For example, we want to construct such linear equation for $\hat U^{of}(\lambda_1;x)$ that after multiplying it by $G^1(\lambda_1;x)$ and integrating over $\lambda_1$  one gets
\begin{eqnarray}\label{DS}
u^{of}_{x_M} - \Delta u^{of} B^2  + u^{of}u_{12}u_{21} B^2 =0,\;\;\Delta = \sum_{i=1}^{M-1} \partial_{x_i}^2,  
\end{eqnarray} 
which becomes multidimensional NLS if $x_M=i t$, $i^2=-1$, $u_{21}=\bar u_{12}$. Let $G^1_{x_1}(\lambda_1;x) = \lambda_1 G^1(\lambda_1;x)$.
Appropriate  linear equation is following
\begin{eqnarray}\label{DS_lin}
\hat U^{of}_{x_M}(\lambda_1;x) - \Delta \hat U^{of}(\lambda_1;x) B^2 + \hat U^{of}(\lambda_1;x) u_{12}u_{21} - B^2 \lambda_1^2 \hat U^{of}(\lambda_1;x)- &&\\\nonumber
2 \Big(\lambda_1  \hat U^{of}_{x_1}(\lambda_1;x)+
 \lambda_1 B^2  \hat U^{of}_{x_2}(\lambda_1;x)+  \lambda_1^2 \hat U^{of}(\lambda_1;x)\Big) B^2&=&0
\end{eqnarray}
Thus nonlinear eq.(\ref{DS}) is equivalent to linear eq.(\ref{DS_lin})
where $\hat U$ is expressed in terms of  the dressing functions by the system (\ref{Sec4:U},\ref{DS_Sec6:Fx}-\ref{DS_c^n}). Detailed discussion of this relation is given in the next subsection.

\subsection{Analysis of the system (\ref{Sec4:U},\ref{DS_Sec6:Fx}-\ref{DS_c^n},\ref{DS_lin} )}
\label{DS_IDE_DS}

In this section we characterize solution space of nonlinear equation
(\ref{DS}) in terms of the dressing functions.  
First step is solving equations
(\ref{DS_Sec6:Fx}-\ref{DS_c_x}) 
for $\Psi(\lambda,\mu;x)$, $\Phi(\lambda;x)$ and $C(\mu;x)$. 
Eqs. (\ref{DS_Sec6:Fx}) represent  nonhomogeneous system  
for $\Psi(\lambda,\mu;x)$, so, similar to Sec.\ref{IDE}, we take the following solution:
\begin{eqnarray}\label{DS_Sec5:F}
\Psi_{\alpha\beta}(\lambda,\mu;x) = \partial_{x_1}^{-1}
\Big(\Phi_\alpha(\lambda;x) 
C_{\alpha\beta}(\mu;x)\Big) +
\delta_{\alpha\beta}\delta(\lambda-\mu) 
e^{-\sum\limits_{j=2}^M \Big(h^j_\alpha(\lambda) +
g^j_{\alpha\beta}(\mu)\Big)},
\end{eqnarray}
Solutions of eqs. (\ref{DS_Phi_x},\ref{DS_c_x})
in view of (\ref{CGG}) read
\begin{eqnarray}\label{DS_Sec4:chi_sol_t}
\label{DS_Phi_def}
\Phi_{\alpha}(\lambda;x)&=& 
e^{K^\Phi_{\alpha}(\lambda;x)},\;\;
K^\Phi_{\alpha}(\lambda,x)=\lambda_2 (x_1 + x_2 B^2_\alpha) + x_M \lambda_2^2 B^2_\alpha  -\sum_{j=2}^{M}
  h^j_\alpha(\lambda)  
 x_j,\\\label{DS_Sec4:c_sol}
 C_{\alpha\beta}(\mu;x)&=& G^1_\alpha(\mu_1;x) G^2_{\alpha\beta}(\mu;x)
 ,\\\nonumber
G^1_\alpha(\mu_1;x) &=&
e^{K^{G^1}_\alpha(\mu_1;x)} ,\;\;
G^2_{\alpha\beta}(\mu;x)= 
 e^{K^{G^2}_{\alpha\beta}(\mu;x)} C^0_{\alpha\beta}(\mu) 
,\;\;\\\nonumber
&& K^{G^1}_\alpha(\mu_1;x)= \mu_1(x_1 + x_2 B^2_\alpha) - x_M \mu_1^2 B^2_\alpha,\;\;
K^{G^2}_{\alpha\beta}(\mu;x)=-\sum_{j=2}^{M}
     g^j_{\alpha\beta}(\mu)  x_{j}.
 \end{eqnarray}
 Thus expression
 (\ref{DS_Sec5:F})  may be written in  explicit
 form:
\begin{eqnarray}\label{DS_Sec6:Fexpl}
\Psi_{\alpha\beta}(\lambda,\mu;x) &=&  
\frac{e^{K^\Phi_{\alpha}(\lambda,\lambda_1;x)+
K^{G^1}_{\alpha}(\mu_1;x)+
K^{G^2}_{\alpha\beta}(\mu;x)} C^0_{\alpha\beta}(\mu) 
 }{\lambda_2+\mu_1} 
 + 
\delta_{\alpha\beta}\delta(\lambda-\mu)
 e^{-\sum\limits_{j=2}^M \Big(h^j_\alpha(\lambda) +
g^j_{\alpha\beta}(\mu)\Big)}.
\end{eqnarray}
Equations (\ref{Sec4:E0^U}-\ref{hatU^2n}) have the same form  
with
\begin{eqnarray}\label{DS_Phi^1} 
\Phi^1_{\alpha\beta}(\lambda;x)&=&
e^{K^{\Phi^1}_{\alpha}(\lambda;x)}  ,
\;\;
K^{\Phi^1}_{\alpha}(\lambda;x)= \lambda_2 (x_1+x_2 B^2) + x_M\lambda_2^2  B^2+   \sum_{j=2}^{M}
 g^j_{\alpha\alpha}(\lambda)  
 x_j.
\end{eqnarray}
Function $\hat \Phi^{of}$  satisfies equation (\ref{DS_lin}) where $\hat U$ is related with $\hat \Phi$ by  eq.(\ref{hatU}).
Note, that diagonal elements of $\hat\Phi$  may be arbitrary functions of single independent  variable, similar to Sec.\ref{DR_nl}.
In particular case $C^0(\mu)=\delta(\mu_1) \tilde C^0(\tilde \mu)$, eq.(\ref{hatU}) 
  reduces to PDE (\ref{RhatU}) so that $u$ is defined by the formula (\ref{potentials}a).
We see that eq.(\ref{DS_lin}) in view of (\ref{RhatU}) is linear  PDE for $\varphi^{of}(\lambda_1;x)$ with "boundary" function $\varphi(0;x)$ satisfying (\ref{potentials}a).
Remark made in the end of Sec.\ref{IDE} regarding numerical construction of 
particular solutions is relevant for this section as well.

\section{Conclusions}
\label{conclusions}
We applied a variant of the dressing method to derive a special representation for  a largely arbitrary multidimensional    
 nonlinear PDEs nonintegrable in classical sence.  
 Although we have considered only $N$-wave equation and NLS, reducible from the linear eqs.(\ref{lin:U}) and (\ref{DS_lin}) respectively, different linear equation for the spectral function $\hat U(\lambda_1;x)$  may be used. The only requirement is that after multiplying this equation by $G^1(\lambda_1;x)$ and integrating over $\lambda_1$ one gets nonlinear PDE for  $u^{of}$. 

We introduced several modifications in the classical dressing
method:
\begin{enumerate}
\item
Eqs.(\ref{Sec6:Fx}) (or(\ref{DS_Sec6:Fx})) with functions 
 $h^n(\lambda)$ and $g^n(\mu)$ showing that derivatives
$\Psi_{x_j}(\lambda,\mu;x)$ are not 
separated functions of spectral parameters.
\item
Eq.(\ref{CGG}) splitting $C(\lambda;x)$. 
\item
Extra constrain
(\ref{condition}) (or (\ref{DS_lin}) together with (\ref{hatU})) defining  structure of
PDE (\ref{nl}) (or (\ref{DS})). This constrain is equation for   function $\hat
\Phi^{of}(\lambda_1,x)$  (see for instance eqs.(\ref{condition_f}, \ref{condition_f2}) of the Sec.\ref{IDE})
and  has   no spectral
origin. 
\end{enumerate}
At the present form, multidimensional dressing method doesn't give explicit 
solutions for nonlinear PDEs, but represents them in different form. We expect perspective development of the ideas outlined in this paper.

Author thanks Prof. P.M.Santini for useful discussions of some
aspects of this paper.
The work was supported by INTAS Young Scientists Fellowship Nr. 04-83-2983,
 RFBR grant 04-01-00508 and grant NSh 1716-2003.

\end{document}